\documentclass[12pt]{article}
\usepackage[latin1]{inputenc}
\usepackage{amsmath, amssymb, eucal, mathrsfs}
\newtheorem{theorem}{Theorem}

\newtheorem{lemma}[theorem]{Lemma}

\newtheorem{corollary}[theorem]{Corollary}

\newtheorem{remark}[theorem]{Remark}

\newtheorem{example}[theorem]{Example}

\newcommand{\R}{\mathbb{R}}
\newcommand{\Les}{\mathbb{L}}
\newcommand{\Q}{\mathbb{Q}}
\newcommand{\Sf}{\mathbb{S}}

\newcommand{\Hy}{\mathbb{H}}
\newcommand{\spa}{\mbox{span\,}}
\newcommand{\hess}{\mbox{Hess\,}}

\newcommand{\kerl}{\mbox{ker}}
\newcommand{\grad}{\mbox{grad\,}}

\newcommand{\po}{{\hspace*{-1ex}}{\bf .  }}

\newcommand{\nap}{\nabla^{\perp}}
\def\<{{\langle}}
\def\>{{\rangle}}

\def\n{\nabla}

\def\a{\alpha}
\def\be{\begin{equation} }
\def\ee{\end{equation} }

\def\nap{\nabla^\perp}
\def\proof{\noindent\emph{Proof: }}
\def\qed{\ifhmode\unskip\nobreak\fi\ifmmode\ifinner
\else\hskip5 pt \fi\fi\hbox{\hskip5 pt \vrule width4 pt
height6 pt  depth1.5 pt \hskip 1pt }}

\newcommand\blfootnote[1]{%
  \begingroup
  \renewcommand\thefootnote{}\footnote{#1}%
  \addtocounter{footnote}{-1}%
  \endgroup
}

\begin{document}

\title{Kaehler submanifolds of hyperbolic space}
\author{M.\ Dajczer and Th.\ Vlachos}
\date{}
\maketitle

\blfootnote{\textup{2010} \textit{Mathematics Subject Classification}:
53B25,  53C21, 53C40.}


\begin{abstract} 
We present several local and global results on isometric immersions of 
Kaehler manifolds $M^{2n}$ into hyperbolic space $\Hy^{2n+p}$. For instance, 
a classification is given in the case of dimension $n\geq 4$ and 
codimension $p=2$.  Moreover, as corollaries of general results, we conclude 
that there are no isometric immersion in codimension $p\leq n-2$ if the Kaehler 
manifold is of dimension $n\geq 4$ and either has a point of positive holomorphic sectional 
curvature or is compact.
\end{abstract}

Since the pioneering work of Dajczer and Gromoll \cite{DG1}, \cite{DG2}, 
\cite{DG3}, \cite{DG4} on real Kaehler submanifolds, that is, isometric 
immersions of Kaehler manifolds into Euclidean space, many authors worked on 
the subject in both the local and global case. For instance, see
\cite{APS}, \cite{BEFT}, \cite{DR}, \cite{Di}, \cite{EFT}, \cite{FZ1}, 
\cite{FZ2}, \cite{FZ3}, \cite{FHZ}, \cite{Fur}, \cite{He}, \cite{MR}, 
\cite{Ry}, \cite{YZ1} and \cite{YZ2}.

A strong result when the ambient space is the round sphere $\Sf^N$ is due 
to Florit, Hui and  Zheng \cite{FHZ}. By taking advantage of the umbilical 
inclusion of the sphere into Euclidean space they proved that any isometric 
immersion $f\colon M^{2n}\to\Sf^{2n+p}$ of a Kaehler manifold with 
codimension $p\leq n-1$ is part of the product of round two-spheres, namely,  
$M^{2n}\subset\Sf^2\times\cdots\times\Sf^2\subset\Sf^{3n-1}\subset\R^{3n}$.

The purpose of this paper is to take aim at the study of isometric immersions 
of Kaehler manifolds  $M^{2n}$, $n\geq 2$, into hyperbolic space $\Hy^{2n+p}$. 
This case is certainly harder than the spherical case, in good  part due to 
the fact the Euclidean space can be isometrically immersed in hyperbolic space 
with codimension one as an umbilical horosphere. Hence, any euclidean submanifold 
becomes an hyperbolic submanifold with codimension one higher. Nevertheless, 
two results have already been obtained in situations that avoid this difficulty. 
In the hypersurface case Ryan \cite{Ry} showed that the only possibility  other 
than the horosphere is $M^4=\Hy^2\times\Sf^2\subset\Hy^5\subset\Les^6$, where 
$\Les^m$  stands for the $m$-dimensional Lorentz space. Dajczer and Rodr\'iguez 
\cite{DR} proved that if we require the immersion to be  minimal then, regardless of 
the codimension, there are no other possibilities  
than minimal surfaces.
\vspace{1ex}

We first consider the local situation in the case of codimension two.

\begin{theorem}\po\label{th1} Let $f\colon M^{2n}\to\Hy^{2n+2}$, $n\geq 4$, 
be an isometric immersion of a Kaehler manifold without flat points. Then 
$f=i\circ g\colon M^{2n}\to\Hy^{2n+2}$ is locally a composition of isometric 
immersions where $g\colon M^{2n}\to\R^{2n+1}$ is a real Kaehler hypersurface 
and $i\colon\R^{2n+1}\to\Hy^{2n+2}$ the inclusion as a horosphere.  
\end{theorem}

It was shown by Dajczer and Gromoll \cite{DG1} that any real Kaehler 
hypersurface without flat points $f\colon M^{2n}\to\R^{2n+1}$, $n\geq 2$, 
can be  locally parametrized by the so called Gauss parametrization in 
terms of a pseudoholomorphic spherical surface $h\colon L^2\to\Sf^{2n}$ 
and a function in $C^\infty(L)$. Calabi \cite{Ca} established a 
correspondence between these surfaces and holomorphic maps into the 
hermitian symmetric space $\wp_n=SO(2n+1)/U(n)$ of all oriented 
hyperplanes in $\R^{2n+1}$ with complex structure. Then Dajczer
and Vlachos \cite{DV} gave a Weierstrass type representation for the
surfaces and showed how this can be used to parametrize the 
hypersurfaces themselves. The trivial case, namely, when $h$ is a 
totally geodesic sphere, corresponds to cylinders where 
$M^{2n}=M^2\times\R^{2n-2}$ and $f=k\times I$ where $k\colon M^2\to\R^3$ 
is any surface and $I$ is the identity map on $\R^{2n-2}$.  These 
submanifolds are the only ones in the class that can be complete 
manifolds.

\begin{example}\po\label{ex}
{\em Theorem \ref{th1} is sharp since it does not hold for $n=3$, as shown by 
$$
M^6=\Hy^2\times\Sf^2\times\Sf^2\subset\Hy^8\subset\Les^9
=\Les^3\times\R^3\times\R^3
$$
where $\Hy^2\subset \Les^3$  and $\Sf^2\subset\R^3$.
}\end{example}

Next we consider the case of submanifolds with higher codimension.
We have the following consequence of a general result given later.

\begin{theorem}\po\label{con}
 If a  Kaehler manifold $M^{2n}$, $n\geq 3$,  has positive
holomorphic sectional curvature at some point then there is 
no isometric immersion in $\Hy^{2n+p}$ for $p\leq n-2$.
\end{theorem}

The Omori-Yau maximum principle for the Hessian is said to
hold on a Riemannian manifold $M^n$ if  for any function 
$g\in C^2(M)$ with $g^*=\sup_Mg<+\infty$ there exists a sequence 
of points $\{x_k\}_{k\in\mathbb N}$ in $M^n$ satisfying:
$$
(i)\; g(x_k)>g^*-1/k,\;\;(ii)\;\|\grad g(x_k)\|<1/k,
\;\;(iii)\;\hess g(x_k)(X,X)\leq (1/k)\|X\|^2
$$ 
for all $X\in T_{x_k}M$. It is well known \cite{amr} that this 
maximum principle holds on a manifold $M^n$ if its sectional curvature 
satisfies
$$
K_M(x)\geq -C\rho^2(x)\left(\Pi_{j=1}^N\log^{(j)}(\rho(x))\right)^2,
\;\;\;\rho(x)>\!\!>1,
$$
for a constant $C>0$, where   $\rho$ is the distance function in 
$M^n$ to a reference point and $\log^{(j)}$ stands for the $j$-iterated logarithm.
\vspace{1ex}

In this paper we  use a weaker version of the above maximum principle.
The \emph{weak maximum principle for the Hessian} amounts to require only 
conditions $(i)$ and $(iii)$.  It is known \cite{Al} that this principle holds 
if $M^n$ is a complete manifold and there  exist a function $\varphi\in C^2(M)$ 
and a constant $k>0$ such that $\varphi(x)\rightarrow +\infty$ as 
$x\rightarrow\infty$ and 
$$
\hess\varphi(\,,\,)\leq k\varphi\<\,,\,\>
$$
outside a compact subset of $M^n$.

It was shown by Mari and Rigoli \cite{MR} that if a Kaehler manifold 
$M^{2n}$ satisfies the weak maximal principle for the Hessian, then it 
cannot be isometrically immersed in a  nondegenerate cone of 
$\R^{3n-1}$. This generalizes the result of Hasanis \cite{Has} who 
assumed completeness and sectional curvature bounded from below to 
conclude that the submanifold must be unbounded.

\begin{theorem}\label{th4}\po Let $f\colon M^{2n}\to\Hy^{2n+p}$, 
$2\leq p\leq n-2$, be an isometric immersion of a Kaehler manifold. 
If the weak principle for the Hessian holds on $M^{2n}$ then $f(M)$ 
is unbounded.
\end{theorem}

In particular, we have the following consequence.

\begin{corollary}\po
There is no isometric immersion of a compact Kaehler manifold
$M^{2n}$ into $\Hy^{2n+p}$ if $n\geq 3$ and $p\leq n-2$.
\end{corollary}

We also consider the local cases of submanifolds of dimensions four 
and six but under the additional assumption of flat normal bundle.

\begin{theorem}\label{fnb}\po Let $f\colon M^6\to\Hy^8$ be an isometric 
immersion with flat normal bundle of a Kaehler manifold without flat 
points. Then $f$ is locally a composition of isometric immersions as 
in Theorem \ref{th1} or is as in Example \ref{ex}.
\end{theorem}

The result for dimension four will be given at the end of next section.

\section{The proofs}

In the sequel, we prove the results stated in the introduction. The 
proof of Theorem \ref{con} is immediate from the following result.

\begin{theorem}\po\label{th2} Let $f\colon M^{2n}\to\Hy^{2n+p}$, $p\leq n-2$, 
be an isometric immersion of a Kaehler manifold. Then at each point
$x\in M^{2n}$ there is a $J$-invariant subspace $L^{2m}\subset T_xM$ 
with $m\geq n-p-1$ such that the  holomorphic sectional curvature  for 
any complex plane $P\subset L^{2m}$ satisfies $K_P\leq 0$. 
\end{theorem}

\proof Recall that the curvature tensor of a Kaehler manifold 
$(M^{2n},J)$ satisfies 
\be\label{kaehler}
R(X,Y)\circ J=J\circ R(X,Y)
\ee
for all $X,Y\in TM$.
Let $\a\colon TM\times TM\to N_fM$ denote the second fundamental 
form of $f$ taking values in the normal bundle.
At $x\in M^{2n}$ let
$$
\beta\colon T_xM\times T_xM\to W=N_fM(x)\oplus N_fM(x)\oplus\R^2
$$ 
be the bilinear form defined by
$$
\beta(X,Y)=(\a(X,Y),\a(X,JY),\<X,Y\>,\<X,JY\>).
$$
If we endow $W$ with the indefinite inner product defined as 
$$
\<\!\<(\xi_1,\eta_1,t_1,s_1),(\xi_2,\eta_2,t_2,s_2)\>\!\>
=\<\xi_1,\xi_2\>-\<\eta_1,\eta_2\>-t_1t_2+s_1s_2,
$$
then the bilinear form $\beta$ is a flat.  This means that
$$
\<\!\<\beta(X,Y),\beta(Z,V)\>\!\>=\<\!\<\beta(X,V),\beta(Z,Y)\>\!\>
$$
for all $X,Y,Z,V\in T_xM$.  In fact, a straightforward computation 
making use of the Gauss equation 
\be\label{gauss}
\<R(X,Y)Z,V\>= -\<(X\wedge Y)Z,V\>+\<\a(X,V),\a(Y,Z)\>
-\<\a(X,Z),\a(Y,V)\>
\ee
yields
$$
\<\!\<\beta(X,Y),\beta(Z,V)\>\!\>-\<\!\<\beta(X,V),\beta(Z,Y)\>\!\>
=\<R(X,Z)V,Y\>-\<R(X,Z)JV,JY\>
$$
where the right hand side vanishes due to \eqref{kaehler}. 

Given $Y\in T_xM$, let $B_Y\colon T_xM\to W$ be defined by
$$
B_Y={\beta}(Y,\,\cdot\,).
$$
A vector $Z\in T_xM$ is said to be (left) regular element of $\beta$ if
$$
\dim B_Z(T_xM)=\max\{\dim B_X(T_xM):X\in T_xM\}.
$$
It has been shown in \cite{Mo} that the subset $RE(\beta)\subset T_xM$ 
of regular elements of a flat bilinear form is open and dense, 
and that
\be\label{moore}
\beta(X,Y)\in B_Z(T_xM)\cap(B_Z(T_xM))^\perp
\ee
for any $X\in T_xM$ and $Y\in\kerl B_Z$.

Fix $Z_0\in RE(\beta)$ and consider the $J$-invariant tangent 
subspace $L=\kerl B_{Z_0}$, that is,
$$
L=\{X\in T_xM:\a(Z_0,X)=\a(Z_0,JX)=\<Z_0,X\>=\<JZ_0,X\>=0\}.
$$
Then 
$$
\dim L\geq 2(n-p-1)\geq 2.
$$
We have from \eqref{moore} that
\be\label{regel}
\<\!\<\beta(X_1,Z_1),\beta(X_2,Z_2)\>\!\>=0
\ee
for all $X_1,X_2\in T_xM$ and $Z_1,Z_2\in L$. 

Let $\gamma\colon T_xM\times T_xM\to N_fM(x)\oplus\R$ be the 
bilinear form given by
$$
\gamma(X,Y)=(\a(X,Y),\<X,Y\>)
$$
where $N_fM(x)\oplus\R$ is endowed with the Lorentzian inner product
$$
\<(\xi_1,t_1),(\xi_2,t_2)\>=\<\xi_1,\xi_2\>-t_1t_2.
$$
Then \eqref{regel} is equivalent to 
\be\label{equiv}
\<\gamma(X_1,Z_1),\gamma(X_2,Z_2)\>=\<\gamma(X_1,JZ_1),\gamma(X_2,JZ_2)\>
\ee
for all $X_1,X_2\in T_xM$ and $Z_1,Z_2\in L$. 

We claim that the inner product induced on the subspace 
$S\subset N_fM(x)\oplus\R$ defined by
$$
S=\spa\{\gamma(X,Z)\;\mbox{for all}\;X\in T_xM,\;Z\in L\}
$$
is degenerate, that is, $S\cap S^\perp\neq 0$. Assume to the 
contrary, and define  $\bar J\colon S\to S$ by
$$
\bar J\gamma(X,Z)=\gamma(X,JZ).
$$
Then $\bar J$ is an isometry from \eqref{equiv} and $\bar J^2=-I$. 
Hence
$$
\gamma(JZ_1,JZ_2)=-\gamma(Z_1,Z_2)
$$
for any $Z_1,Z_2\in L$.  In particular, this gives
$$
\<JZ_1,JZ_2\>=-\<Z_1,Z_2\>,
$$
which is a contradiction and proves the claim.

That $S$ is degenerate means that there is a  unique unit vector 
$\delta\in N_fM(x)$ such that $(\delta,1)\in S$ and
\be\label{deg}
S\cap S^\perp=\spa\{(\delta,1)\}.
\ee
Let $\phi\colon T_xM\times T_xM\to R=(\spa\{\delta\})^\perp\subset N_fM$
be defined by
\be\label{sff}
\a(X,Y)=\phi(X,Y)+\<A_\delta X,Y\>\delta
\ee
where $A_\delta$ is the shape operator associated to $\delta$. Hence
$$
\gamma(X,Y)=\phi(X,Y)+\<A_\delta X,Y\>\delta+\<X,Y\>(0,1)
$$
for all $X,Y\in T_xM$. 

We have that \eqref{deg} is equivalent to 
\be\label{deg2}
\<A_\delta X,Z\>=\<X,Z\>
\ee
for any $X\in T_xM$ and $Z\in L$. Then \eqref{equiv} is equivalent to
\be\label{newj}
\<\phi(X_1,Z_1),\phi(X_2,Z_2)\>=\<\phi(X_1,JZ_1),\phi(X_2,JZ_2)\>
\ee
for any $X_1,X_2\in T_xM$ and $Z_1,Z_2\in L$.  Set 
$$
R_0=\spa\{\phi(X,Z)\;\mbox{for all}\;X\in T_xM,\;Z\in L\}.
$$
Then the map $\tilde J\colon R_0\to R_0$ defined by
$$
\tilde J\phi(X,Z)=\phi(X,JZ)
$$
is an isometry and satisfies $\tilde J^2=-I$. In particular,
\be\label{phi}
\phi(JZ_1,JZ_2)=-\phi(Z_1,Z_2)
\ee
for any $Z_1,Z_2\in L$. 

It is now easily seen using \eqref{sff}, \eqref{deg2} and \eqref{phi} that the Gauss 
equation \eqref{gauss} gives for the sectional curvature of the 
holomorphic plane $P=\spa\{Z,JZ\}$ that
$$
K_P=-\|\phi(Z,Z)\|^2-\|\phi(Z,JZ)\|^2\leq 0,
$$
for any $Z\in L$, and this concludes the proof.\vspace{2ex}\qed

\noindent\emph{Proof of Theorem \ref{th1}:}
We use the notation and arguments from the proof of Theorem \ref{th2}. 
Since the codimension is $p=2$, then  $\dim L\geq 2n-6\geq 2$. From 
\eqref{newj} we have
$$
\<\phi(X,Z),\phi(X,JZ)\>=0
$$
for any $X\in T_xM$ and $Z\in L$. Hence $\phi(X,Z)=0$ by dimension reasons. 
Thus
$$
\a(X,Z)=\<X,Z\>\delta
$$
for any $X\in T_xM$ and $Z\in L$. Now the Gauss equation gives
$$
\<R(X,Z)Z,X\>=-1+\<\a(X,X),\a(Z,Z)\>-\|\a(X,Z)\|^2=-1+\<A_\delta X,X\>
$$
and 
$$
\<R(X,Z)JZ,JX\>=\<\a(X,JX),\a(Z,JZ)\>-\<\a(X,JZ),\a(Z,JX)\>=0
$$
for any unit vectors $X\in L^\perp$ and $Z\in L$. We conclude
using \eqref{kaehler} that $A_\delta=I$.

We obtain from \eqref{deg} that $\delta$ can be taken locally 
to be a smooth vector field. If this umbilical unit vector field 
is parallel in the normal connection on an open connected subset 
$U\subset M^{2n}$, then an elementary argument gives that $f(U)$ 
is contained in an horosphere.  
If $\delta$ is not parallel at some point of 
$M^{2n}$, then the same holds in some open neighborhood 
$V\subset M^{2n}$. By the Codazzi equation, we have
$$
A_{\nap_X\delta}Y=A_{\nap_Y\delta}X.
$$
It follows easily that $\dim\kerl A_{\delta^\perp}\geq 2n-1$.  Hence, from
the Gauss equation we have that $V$ is flat, and this is a contradiction.
\qed

\begin{remark}\po{\em  The existence of the umbilical vector field
$\delta$ in the proof above can also be obtained from Lemma $3$ in 
\cite{DG4}}.
\end{remark}

\noindent\emph{Proof of Theorem \ref{th4}:}
Let $r\colon\Hy^{2n+p}\to[0,+\infty)$ be the distance function
in $\Hy^{2n+p}$ to a reference point. Its gradient satisfies 
\be\label{norm}
\|\grad^H r\|=1.
\ee 
Moreover, the Hessian comparison theorem yields
\be\label{compar}
\hess^H r(Y,Y)(x)=\coth(r(x))\left(\|Y\|^2-\<\grad^H r,Y\>^2\right)
\ee
for any $Y\in T_x\Hy^{2n+p}$.

Given $h\in C^\infty(\Hy^{2n+p})$, it is a standard fact that 
then the Hessians of $h$ and $g=h\circ f$ are related by
\be\label{rel}
\hess g(x)(X,Y)=\hess^H h(f(x))(f_*X,f_*Y)+\<\grad^H h(f(x)),\a(X,Y)\>
\ee
for any $X,Y\in T_xM$.
\vspace{1ex}

Set $g=\psi\circ r\circ f$ where $\psi\in C^\infty(\R)$. Thus
$$
\hess g(x)(X,X)=\hess^H \psi\circ r(f(x))(f_*X,f_*X)
+\psi'(r(f(x)))\<\grad^Hr(f(x)),\a(X,X)\>
$$
for any $X\in T_xM$. On the other hand, 
\begin{align*}
\hess^H\psi\circ r(f(x))(f_*X,f_*X)
=&\,\psi^{''}(r(f(x)))\<\grad^Hr(f(x)),f_*X\>^2\\
&+\psi'(r(f(x)))\hess^Hr(f(x))(f_*X,f_*X)
\end{align*}
for any $X\in T_xM$. Setting $\psi(t)=\cosh t$, it follows using \eqref{compar} 
and \eqref{rel} that
\be\label{final}
\hess g(x)(X,X)
=\cosh(r(f(x)))\|X\|^2+\sinh(r(f(x)))\<\grad^H r(f(x)),\a(X,X)\>
\ee
for any $X\in T_xM$.

We obtain from \eqref{sff} and \eqref{deg2} that
$$
\<\grad^H r(f(x)),\a(Z,Z)\>=\<\grad^H r(f(x)),\phi_x(Z,Z)\>
+\<\grad^H r(f(x)),\delta_x\>\|Z\|^2
$$
for all $Z\in L(x)$. Moreover, we have from \eqref{phi} that 
$$
\phi_x(Z,Z)+\phi_x(JZ,JZ)=0
$$
for all $Z\in L(x)$. Hence \eqref{final} gives
\begin{align}\label{*}
\hess g(x)(Z,Z)&+\hess g(x)(JZ,JZ)\nonumber\\
&=2\left(\cosh(r(f(x)))
+\sinh(r(f(x)))\<\grad^H r(f(x)),\delta_x\>\right)\|Z\|^2
\end{align}
for all $Z\in L(x)$.

Suppose that $f(M)$ is bounded. By assumption there exists a 
sequence of points $\{x_k\}_{k\in\mathbb N}$ in $M^{2n}$ such that
$$
\lim g(x_k)=g^*=\sup g<+\infty\;\;\mbox{and}\;\;
\hess g(x_k)(X_k,X_k)\leq\frac{1}{k}\|X_k\|^2, X_k\in T_{x_k}M.
$$
Setting $r_k=r(f(x_k))$ and letting $k\to+\infty$,  we have that 
$g^*=\lim\cosh(r_k)$. Therefore $\lim r_k=r^*=\sup\{r\circ f\}>0$. 
We obtain from \eqref{*}  that
$$
\frac{1}{k}\geq\cosh r_k+\sinh r_k\<\grad^H r(f(x_k)),\delta_{x_k}\>.
$$
Taking the limit as $k\to+\infty$ and using \eqref{norm} gives
$$
0\geq \cosh r^*-\sinh r^*>0,
$$
and that is a contradiction. \vspace{1ex}\qed

For the proof of Theorem \ref{fnb} we need the following lemmas.

\begin{lemma}\po\label{var}  Let $f\colon M^{2n}\to\Q_c^{2n+p}$ be an 
isometric immersion of a Kaehler manifold  with flat normal bundle into 
a space form of sectional curvature $c$. Then
\be\label{uno}
\<JX_i,X_k\>K_{ij}=0,\;\;i\neq j\neq k\neq i,
\ee
where $X_1,\ldots,X_{2n}$ is an orthonormal tangent frame that diagonalizes 
$\a$ and $K_{ij}$ denotes the sectional curvature of the plane 
$\spa\{X_i,X_j\}$. 
Moreover, if $K_{ij}\neq 0$ then:
\begin{itemize}
\item[(i)] The plane $\spa\{X_i,X_j\}$ is holomorphic.
\item[(ii)] The sectional curvatures satisfy
\be\label{i}
K_{ik}=K_{jk}=0\;\;\mbox{for all}\;\;k\not\in\{i,j\}.
\ee
\item[(iii)] The second fundamental form satisfies
\be\label{ii}
\<\a_i-\a_j,\a_k\>=0\;\;\mbox{for all}\;\;k\not\in\{i,j\}
\ee
where $\a_j=\a(X_j,X_j)$. 
\end{itemize}
\end{lemma}

\proof  We obtain from \eqref{kaehler} that
\be\label{curv}
\<R(X_i,X_j)JX_i,X_k\>+\<R(X_i,X_j)X_i,JX_k\>=0
\;\;\mbox{if}\;\;k\not\in\{i,j\}.
\ee
Computing each term of \eqref{curv} by means of the Gauss 
equation gives that \eqref{uno} and \eqref{curv} are equivalent. 
Parts $(i)$ and $(ii)$ follow from \eqref{uno}.  
Then
$$
0=K_{ik}=c+\<\a_i,\a_k\>,\;\;0=K_{jk}=c+\<\a_j,\a_k\>,
$$
and also part  $(iii)$ follows.\qed

\begin{lemma}\po With the above notations
the  Codazzi equations of $f$ are equivalent to
\be\label{cod1}
\nap_{X_i}\a_j=\Gamma_{jj}^i(\a_j-\a_i),\;\;j\neq i,
\ee
and
\be\label{cod2}
\Gamma_{ij}^k(\a_k-\a_j)=\Gamma_{ji}^k(\a_k-\a_i),\;\;k\neq i\neq j\neq k,
\ee
where $\Gamma_{ij}^k=\<\n_{X_i}X_j,X_k\>$. 
\end{lemma}

\proof A straightforward computation.\vspace{2ex}\qed

\noindent\emph{Proof of Theorem \ref{fnb}:}
To conclude that $f$ is as in Theorem \ref{th1}, it is well-known 
that one has to show that there is a unit normal vector field $\xi$ 
that is  parallel in the normal connection and such that $A_\xi=I$.
\vspace{1ex}

We will consider two cases:
\vspace{1ex}

\noindent Case $(I)$.  Assume $K_{ij}\neq 0$ and $\a_i-\a_j\neq 0$ 
for some $i\neq j$. We claim that
$$
K_{rs}=0\;\;\mbox{and}\;\;\a_r=\a_s=\xi\perp (\a_i-\a_j)
\;\;\mbox{with}\;\;\|\xi\|=1\;\;\mbox{for}\;\;r,s\not\in\{i,j\}.
$$
Say $i=1$, $j=2$. From \eqref{ii} we have $\a_r=\mu\a_s$. Then \eqref{i} 
and the Gauss equation yield
$$
\<\a_k,\a_s\>=1=\<\a_k,\a_r\>,\;\;k=1,2.
$$
In particular $\a_1,\a_2$ are linearly independent.
Hence $\mu=1$ and $\a_r=\xi$ if $r\not\neq\{1,2\}$. Therefore
$$
K_{rs}=-1+\|\xi\|^2\;\;\mbox{if}\;\; 3\leq r\neq s\leq 6. 
$$
We now obtain from part $(ii)$ of Lemma \ref{var} that $K_{rs}=0$, and 
the claim follows. 

We have 
$$
\<\a_i,\xi\>=1\;\;\mbox{if}\;\; i=1,2.
$$
Hence $\xi$ is an umbilical vector field. Moreover, using 
$K_{12}\neq 0$ we obtain that $\xi\neq\a_1,\a_2$. It now follows from 
\eqref{cod2} that
$$
\Gamma_{ri}^s=0\;\;\mbox{if}\;\;r\neq s\geq 3
\;\;\mbox{and}\;\;i=1,2.
$$
Set
$$
D_1=\spa\{X_1,X_2\}\;\;\mbox{and}\;\;D_2=\spa\{X_3,\ldots,X_6\}.
$$
Since $D_1$ is $J$-invariant, then also is $D_2$. Hence
$$
\<\n_{X_r}X_r,X_i\>=\<J\n_{X_r}X_r,JX_i\>=\<\n_{X_r}JX_r,JX_i\>=0\;\;
\mbox{if}\;\; r\geq 3\;\;\mbox{and}\;\; i=1,2.
$$
Thus $D_2$ is totally geodesic distribution. It now follows easily 
from \eqref{cod1}  that $\xi$ is parallel in the normal connection.
Hence $f$ is a composition of isometric immersions as in Theorem \ref{th1}.
\bigskip

\noindent Case $(II)$.  Assume that $K_{ij}\neq 0$ for some $i\neq j$,
and $\a_i=\a_j$. In the sequel, assume $K_{12}\neq 0$, hence $\a_1=\a_2$.
\vspace{1ex}

\noindent Subcase $(a)$.  Assume $K_{rs}=0$ for all $r\neq s\geq 3$.  
The Gauss equation yields
$$
\<\a_1,\a_r\>=1\;\;\mbox{and}\;\;\<\a_r,\a_s\>=1.
$$
It is not difficult to obtain from this set of equations that
$$
\a_r=\xi\neq\a_1\;\;\mbox{with}\;\;\|\xi\|=1.
$$
Thus $\xi$ is an umbilical vector field as in the previous case. 
Moreover, the same type of argument gives that $\xi$ is parallel in 
the normal connection. Again $f$ is a composition as in Theorem \ref{th1}.
\vspace{1ex}

\noindent Subcase $(b)$.  Assume $K_{34}\neq 0$ but $K_{56}=0$.  
Hence
$$
\a_1=\a_2,\;\;\a_3=\a_4,\;\;\<\a_5,\a_6\>=1
$$
and $\a_1-\a_3\neq 0$ since
\be\label{arg}
0\neq K_{12}=-1+\<\a_1,\a_1\>\;\;\mbox{and}
\;\; 0=K_{13}= -1+\<\a_1,\a_3\>.
\ee
 From
$$
\<\a_1,\a_5\>=1=\<\a_1,\a_6\>,\;\;
\<\a_3,\a_5\>=1=\<\a_3,\a_6\>
$$
and $K_{56}=0$ it follows easily that
$$
\a_5=\a_6=\delta\;\;\mbox{with}\;\;\|\delta\|=1.
$$
By an argument as in \eqref{arg} it follows that $\delta\neq\a_1,\a_3$. 
Since 
$$
\<\a_1,\delta\>=1=\<\a_3,\delta\>,
$$
then
$$
\<\a_1-\delta,\delta\>=0=\<\a_3-\delta,\delta\>
$$
Hence
$$
\a_3-\delta=\mu(\a_1-\delta)\;\;\mbox{where}\;\;\mu\neq 0,1.
$$
Since taking the inner product with $\a_1$ yields a contradiction,
this case is not possible.
\vspace{1ex}

\noindent Subcase $(c)$.  Assume $K_{34}\neq 0$ and $K_{56}\neq0$.  
Hence $\a_1=\a_2$, $\a_3=\a_4$ and $\a_5=\a_6$.  Moreover, we have
$\a_1\neq\a_3\neq\a_5\neq\a_1$. In fact, if  $\a_1=\a_3$
then $0=K_{13}=-1+\|\a_1\|^2$ which gives $K_{12}=0$, a contradiction.

As above, we obtain that
\be\label{perp}
\a_i\perp(\a_j-\a_k)\;\;\mbox{if}\;\;i\neq j\neq k\neq i\;\;\mbox{and} 
\;\;i,j,k=1,3,5.
\ee

Similar arguments as above yield that 
$D_1=\spa\{X_1,X_2\}$, $D_2=\spa\{X_3,X_4\}$ and $D_3=\spa\{X_5,X_6\}$
are totally geodesic distributions. From \eqref{cod2} and \eqref{perp} 
we have
$$
\<\n_{D_i}D_j,D_k\>=0\;\;\mbox{if}\;\;i\neq j\neq k\neq i.
$$
Since the distributions are totally geodesic we obtain that
they are parallel, that is,
$$
\n_{D_i}D_j\subset D_j\;\;\mbox{for all}\;\; i,j,k.
$$
It follows that the manifold is locally a product of three surfaces
$M_1^2\times M_2^2\times M_3^2$.  Moreover, since the second fundamental
form of $f$ satisfies $\alpha(D_i,D_j)=0$ if $i\neq j$, then by a result 
due to N\"{o}lker \cite{No} the immersion is an exterior product of 
immersions in $\Hy^8$ as in the statement.\qed

\begin{theorem}\po Let $f\colon M^4\to\Q_c^6$, $c\neq 0$, be an isometric 
immersion with flat normal bundle of a Kaehler manifold free of flat points. 
Then one of the following holds:
\begin{itemize}
\item[(i)] $c=1$ and $f$ is the following external product of immersions.
$$
f=h\times id\colon L^2\times\Sf^2_{c_2}\to\Sf^3_{c_1}\times\Sf^2_{c_2}
\subset\Sf^6_1\subset\R^4\times\R^3
$$
where  $c_1^2+c_2^2=1$.
\item[(ii)] $c=-1$ and $M^4$ is free of flat points. Then either $f$ is a
composition of immersions as in Theorem \ref{th1} or is one of the following 
external product of immersions:
\begin{itemize}
\item[(a)] $f=h\times id\colon L^2\times\Hy^2_{c_1}\to\Sf^3_{c_1}\times\Hy^2_{c_2}
\subset\Hy^6_{-1}\subset\R^4\times\Les^3,\;\;c_1^2-c_2^2=-1$,
\item[(b)] $f=h\times id\colon L^2\times\Sf^2_{c_1}\to\Hy^3_{c_1}\times\Sf^2_{c_2}
\subset\Hy^6_{-1}\subset\Les^4\times\R^3,\;\;c_1^2-c_2^2=1$.
\end{itemize}
\end{itemize}
\end{theorem}

\proof It is omitted since it is quite similar to the one of Theorem \ref{fnb}.

{\renewcommand{\baselinestretch}{1}
\hspace*{-20ex}\begin{tabbing} \indent\= IMPA -- Estrada Dona Castorina, 110
\indent\indent\= Univ. of Ioannina -- Math. Dept. \\
\> 22460-320 -- Rio de Janeiro -- Brazil  \>
45110 Ioannina -- Greece \\
\> E-mail: marcos@impa.br \> E-mail: tvlachos@uoi.gr
\end{tabbing}}
\end{document}